\documentclass{amsart}
\usepackage{amsfonts,amssymb,amsmath,amsthm}
\usepackage{url}
\usepackage{enumerate}

\makeatletter
\@namedef{subjclassname@2010}{  \textup{2010} Mathematics Subject Classification.}
\makeatother
\newtheorem{thrm}{Theorem}[section]

\theoremstyle{definition}
\newtheorem{definition}[thrm]{Definition}
\newtheorem{remark}[thrm]{Remark}
\numberwithin{equation}{section}
\email{hzoubeir2014@gmail.com}
\input{tcilatex}

\begin{document}
\address{ }
\author{Hicham Zoubeir}
\address{Ibn Tofail University, Department of Mathematics,\\
Faculty of Sciences, P. O. B : $133,$ Kenitra, Morocco.}
\title[A note on the convergence of the Adomian decomposition method]{A note
on the convergence of the Adomian decomposition method}

\begin{abstract}
In this note we obtain a new convergence result for the Adomian
decomposition method.
\end{abstract}

\dedicatory{$\emph{This}$ $\emph{modest}$ $\emph{work}$ $\emph{is}$ $\emph{%
dedicated}$ $\emph{to}$ $\emph{the}$ $\emph{memory}$ $\emph{of}$ $\emph{our}$
$\emph{beloved}$ $\emph{master}$ $\emph{Ahmed}$ $\emph{Intissar}$ $\emph{%
(1951-2017),}$ $\emph{a}$ $\emph{distinguished}$ $\emph{professor,}$ $\emph{a%
}$ $\emph{brilliant}$ $\emph{mathematician,}$ $\emph{a}$ $\emph{man}$ $\emph{%
with}$ $\emph{a}$ $\emph{golden}$ $\emph{heart.}$\emph{\ }}
\subjclass[2010]{ 97N40}
\keywords{Adomian decomposition method, Analytic operators.}
\maketitle

\section{ Introduction}

The Adomian decomposition method (ADM) was developped in the $1980$'s by the
American physicist G. Adomian $(1923-1996)$ (\cite{ADM1}-\cite{ADM4}) as a
powerful method for solving functional equations. The main idea of this
method lies in the decomposition of the solution $u$ of a vector nonlinear
equation $u=\mathfrak{f+}N(u)$ where $N$ is an analytic operator and $%
\mathfrak{f}$ a given vector, into a series $u=\underset{n=0}{\overset{%
+\infty }{\sum }}u_{n}$ such that the term $u_{n+1}$ is determined from the
terms $u_{0},u_{1},...,u_{n}$ by a reccurence relation involving a
polynomial $A_{n}$ generated by the Taylor expansion of the operator $N.$
Many papers on the applications of the ADM to the problems arising from
different areas of pure and applied sciences have been published (\cite{AGO1}%
-\cite{AGO3}, \cite{BHAT}, \cite{DUA1}, \cite{DUA2}, \cite{YOUS}). Many
works have been also devoted to the convergence of the ADM (\cite{BABO}, 
\cite{BOUM}-\cite{CHR4}). However let us pointwise that some recent works on
the application of the ADM to nonlinear problems avoid the theoretical
treatment of the question of the convergence of the method. On the other
hand we observe that the current convergence criteria of the ADM give rise
to small regions of convergence (that is to small set of vectors $\mathfrak{f%
}$ for which the vector series $\sum u_{n}$ is convergent according to the
criteria of convergence) which limits the application of this method. Our
purpose in this note is to obtain a new convergence result for the ADM,
which we believe will contribute to strenghting the method by enlarging the
field of its applications.

\section{Preliminary notes}

\subsection{Main definitions}

Along this paper $\left( X,\text{ }\left \Vert .\right \Vert \right) $ is a
given real or complex Banach space.

The following definitions were introduced in (\cite{WHIT}).

\begin{definition}
$\mathfrak{L}_{n}(X)$\ denotes for each $n\in 
\mathbb{N}
^{\ast }$\ a set of continuous symetric $n$-multilinear mappings from $X^{n}$%
\ to $X$. It is well known that $\mathfrak{L}_{n}(X)$\ is a Banach space for
the norm $||.||^{[n]}$\ defined by the relation $:$\ 
\begin{equation*}
||f||^{[n]}:=\sup_{(x_{1},\dots ,x_{n})\in (X-\{0\})^{n}}\left( \frac{%
||f(x_{1},\dots ,x_{n})||}{||x_{1}||\dots ||x_{n}||}\right)
\end{equation*}%
Let $x\in X$\ and $F\in \mathfrak{L}_{n}(X)$. We will denote $F(x,\dots ,x)$%
\ by $F.x^{[n]}$. Given $x_{1},...,x_{p}\in X,$ $n_{1},...,n_{p}\in 
\mathbb{N}
$ such that $n_{1}+...+n_{p}=n,$ we will write $F.\left(
x_{1}^{[n_{1}]}\cdot \cdot \cdot x_{p}^{[n_{p}]}\right) $\ for $F\left( 
\underset{n_{1}\text{times}}{\underbrace{x_{1},...,x_{1}}},\dots ,\underset{%
n_{p}\text{ times}}{\underbrace{x_{p},...,x_{p}}}\right) .$ For convenience
we set $\mathfrak{L}_{0}(X):=X$ and $||.||^{[0]}:=||.||.$
\end{definition}

\begin{definition}
A power series in $x$\ with values in $X$\ (\cite{WHIT}) is a series of the
form $\sum_{n=0}^{\infty }F_{n}.x^{[n]}$, where $F_{0}\in X$\ and $F_{n}\in 
\mathfrak{L}_{n}(X)$\ for $n\geq 1$.
\end{definition}

\begin{definition}
A mapping $F:X\rightarrow X$\ is called an analytic operator on $X$\ (\cite%
{WHIT}) if there exists for each vector $a\in X$\ a sequence $A_{n}(a)\in 
\mathfrak{L}_{n}(X)~~(n\geq 1)$\ and a open neighborhood $U_{a}$\ such that
the real series $\sum \frac{1}{n!}\Vert A_{n}(a)\Vert ^{\lbrack n]}\Vert
x-a\Vert ^{n}$\ is convergent when $x\in U_{a}$\ and the power series $\sum 
\frac{1}{n!}A_{n}(a)(x-a)^{[n]}$\ is convergent to $F(x)-F(a)$\ for all $%
x\in U_{a}$.
\end{definition}

\begin{definition}
The radius of the power series $\sum \frac{1}{n!}A_{n}(a)(x-a)^{[n]}$\ is
the radius of convergence of the complex variable power series\textit{\ }$%
\sum \frac{1}{n!}\left \Vert \mathfrak{A}_{n}(a)\right \Vert ^{\left[ n%
\right] }t^{n}$, that is the number \textit{:}%
\begin{equation*}
\mathfrak{R}\left( F,a\right) :=\frac{1}{\underset{n\rightarrow +\infty }{%
\lim \sup }\left( \frac{1}{n!}\left \Vert \mathfrak{A}_{n}(a)\right \Vert ^{%
\left[ n\right] }\right) ^{\frac{1}{n}}}
\end{equation*}%
with the conventions that $\frac{1}{0}=+\infty $.
\end{definition}

\begin{definition}
Let $F:X\rightarrow X$\ be an analytic operator on the Banach space $X$\ and 
$a\in X.$\ We say that $F$\ has an infinite radius of convergence at the
point $a$\ if $\mathfrak{R}\left( F,a\right) =+\infty $.
\end{definition}

\begin{remark}
If a mapping $F:X\rightarrow X$\ is an analytic operator on $X$\ then the
operator $F$\ is of class $C^{\infty }$\ on $X$\ and the following relations
hold for every $a\in X$\ and $n\in 
\mathbb{N}
:$%
\begin{equation*}
F^{\left( n\right) }\left( a\right) =\frac{1}{n!}\mathfrak{A}_{n}(a)
\end{equation*}
\end{remark}

\begin{remark}
If an analytic operator $F$\ has an infinite radius of convergence at some
point $a\in X$\ on the Banach space $X$, then thanks to $($\cite{WHIT}$)$\ $%
F $\ will have the same property at every point of $X$. Thence we will say,
without loss of precision, that $F$\ has an infinite radius of convergence.
\end{remark}

\section{Abstract presentation of the ADM}

Let $\mathfrak{f}\in X$ and $\mathfrak{N}:X\rightarrow X$ \ an analytic
operator on $X$ with infinite radius of convergence. We consider the vector
equation : 
\begin{equation}
\mathfrak{u}=\mathfrak{N}(\mathfrak{u})+\mathfrak{f}  \label{E_CANON}
\end{equation}%
The ADM for solving the vector equation (\ref{E_CANON}) consists in writing
the unknown vector $\mathfrak{u}$ in the form of an absolutely convergent
vector series $\mathfrak{u}=\sum_{n=0}^{+\infty }\mathfrak{u}_{n}$ and in
splitting the nonlinear term $\mathfrak{N}(\mathfrak{u})$ into an absolutely
convergent vector series $N(\mathfrak{u})=\sum_{n=0}^{+\infty }A_{n}$ where
the term $A_{n}$ is obtained for all $n\in 
\mathbb{N}
\ $by the formula :%
\begin{equation*}
A_{n}:=\frac{1}{n!}\frac{d^{n}}{d\varepsilon ^{n}}\left[ \mathfrak{N}\left( 
\underset{j=0}{\overset{+\infty }{\sum }}\varepsilon ^{j}\mathfrak{u}%
_{j}\right) \right] _{\varepsilon =0}
\end{equation*}%
Thence the equation (\ref{E_CANON}) becomes :%
\begin{equation*}
\underset{n=0}{\overset{+\infty }{\sum }}\mathfrak{u}_{n}=\mathfrak{f}+%
\underset{n=0}{\overset{+\infty }{\sum }}A_{n}
\end{equation*}%
Then we set by formal identification :%
\begin{equation}
\left \{ 
\begin{array}{c}
\mathfrak{u}_{0}=\mathfrak{f} \\ 
\mathfrak{u}_{n+1}=A_{n}%
\end{array}%
\right.  \notag
\end{equation}%
and the following relation holds for every $n\in 
\mathbb{N}
$ :%
\begin{equation}
A_{n}:=\frac{1}{n!}\frac{d^{n}}{d\varepsilon ^{n}}\left[ \mathfrak{N}\left( 
\underset{j=0}{\overset{n}{\sum }}\varepsilon ^{j}\mathfrak{u}_{j}\right) %
\right] _{\varepsilon =0}  \label{RECC-SIMPL}
\end{equation}%
We denote by $\mathcal{A}dm\left( \mathfrak{N};\mathfrak{f}\right) $ the
vector series $\sum \mathfrak{u}_{n}$ if it is well-defined. In the sequel
we will continue to denote by $\mathfrak{u}_{n}$ the general term of the
vector series $\mathcal{A}dm\left( \mathfrak{N};\mathfrak{f}\right) .$

\begin{theorem}
(\cite{ABBA})
\end{theorem} \textit{Let }$\mathfrak{f}\in X,$ $\mathfrak{N:}$ $%
X\rightarrow X$ \textit{be an analytic mapping on }$X$\textit{\ which has an
infinite radius of convergence. Then the following formula holds for every }$%
n\in 
\mathbb{N}
:$%
\begin{equation}
A_{n}=\underset{j_{1}+2j_{2}+...+nj_{n}=n}{\sum }\frac{1}{j_{1}!...j_{n}!}%
\mathfrak{N}^{\left( j_{1}+...+j_{n}\right) }(\mathfrak{f}).\left( \mathfrak{%
u}_{1}^{\left[ j_{1}\right] }\cdot \cdot \cdot \mathfrak{u}_{n}^{\left[ j_{n}%
\right] }\right)  \label{fundamental formula}
\end{equation}

\begin{proof}
Since $\mathfrak{N}$ is analytic, we can write for all $n\in 
\mathbb{N}
$ and $\varepsilon >0:$%
\begin{eqnarray*}
&&\mathfrak{N}\left( \underset{j=0}{\overset{n}{\sum }}\varepsilon ^{j}%
\mathfrak{u}_{j}\right) =\mathfrak{N}(u_{0})+\underset{p=0}{\overset{+\infty 
}{\sum }}\frac{1}{p!}\mathfrak{N}^{\left( p\right) }(\mathfrak{u}%
_{0}).\left( \underset{j=0}{\overset{n}{\sum }}\varepsilon ^{j}\mathfrak{u}%
_{j}\right) ^{\left[ p\right] } \\
&=&\mathfrak{N}(\mathfrak{f})+ \\
&&+\underset{p=0}{\overset{+\infty }{\sum }}\frac{1}{p!}\left( \underset{%
j_{0}+j_{1}+...+j_{n}=p}{\sum }\frac{p!}{j_{1}!...j_{n}!}\varepsilon
^{j_{1}+2j_{2}+...+nj_{n}}\mathfrak{N}^{\left( p\right) }(\mathfrak{f}%
).\left( \mathfrak{u}_{1}^{\left[ j_{1}\right] }\cdot \cdot \cdot \mathfrak{u%
}_{n}^{\left[ j_{n}\right] }\right) \right) \\
&=&\mathfrak{N}(\mathfrak{f})+ \\
&&+\underset{q=0}{\overset{+\infty }{\sum }}\varepsilon ^{q}\left( \underset{%
j_{1}+2j_{2}+...+nj_{n}=q}{\sum }\frac{1}{j_{1}!...j_{n}!}\mathfrak{N}%
^{\left( j_{0}+j_{1}+...+j_{n}\right) }(\mathfrak{u}_{0}).\left( \mathfrak{u}%
_{1}^{\left[ j_{1}\right] }\cdot \cdot \cdot \mathfrak{u}_{n}^{\left[ j_{n}%
\right] }\right) \right)
\end{eqnarray*}%
It follows then from the relation (\ref{RECC-SIMPL}) that the following
relation holds for each $n\in 
\mathbb{N}
$ :%
\begin{equation*}
A_{n}=\underset{j_{1}+2j_{2}+...+nj_{n}=n}{\sum }\frac{1}{j_{1}!...j_{n}!}%
\mathfrak{N}^{\left( j_{0}+j_{1}+...+j_{n}\right) }(\mathfrak{f}).\left( 
\mathfrak{u}_{1}^{\left[ j_{1}\right] }\cdot \cdot \cdot \mathfrak{u}_{n}^{%
\left[ j_{n}\right] }\right)
\end{equation*}
\end{proof}

\section{Statement of the main result}

Our main result in this paper is the following.

\begin{theorem}
\end{theorem}\textit{Let }$\mathfrak{f}\in X,$ $\mathfrak{N:}$ $X\rightarrow
X$ \textit{be an analytic mapping on }$X$\textit{\ on some open neighborhood
of the vector }$\mathfrak{f}$\textit{\ and such that} \textit{the following
estimates hold }$:$%
\begin{equation}
||\mathfrak{N}^{\left( n\right) }(\mathfrak{f})||^{\left[ n\right] }\leq
Ma^{n}n!,\text{ }n\in 
\mathbb{N}
\label{ASSUMP1}
\end{equation}%
\textit{where the constants }$M,$ $a>0$\textit{\ satisfy the condition }:%
\begin{equation}
Ma\leq \frac{1}{4}  \label{ASSUMP2}
\end{equation}

$1.$\textit{Then the vector series }$\mathcal{A}dm\left( \mathfrak{N};%
\mathfrak{f}\right) $\textit{\ is absolutely convergent in the \ Banach
space }$X$ \textit{to} \textit{a vector $\mathfrak{u}$} \textit{which is a
solution of the equation (\ref{E_CANON}). }

$2.$ \textit{If }$Ma=\frac{1}{4},$\textit{then the vector $\mathfrak{u}$
fullfiles the following estimates }$:$%
\begin{equation*}
\left \{ 
\begin{array}{c}
\left \Vert \mathfrak{u}\right \Vert \leq \left \Vert \mathfrak{f}\right
\Vert +\left( 1+\frac{8}{\sqrt{3\pi }}\right) M \\ 
\left \Vert \mathfrak{u-}\underset{j=0}{\overset{n}{\sum }}\mathfrak{u}%
_{j}\right \Vert \leq \frac{8M}{\sqrt{3\pi }}\frac{1}{\sqrt{n}},\mathit{\ }%
n\in 
\mathbb{N}
^{\ast }%
\end{array}%
\right.
\end{equation*}

$3.$ \textit{If }$Ma<\frac{1}{4}$ \textit{then the vector $\mathfrak{u}$
fullfiles the following estimates }$:$%
\begin{equation*}
\left \{ 
\begin{array}{c}
\left \Vert \mathfrak{u}\right \Vert \leq \left \Vert \mathfrak{f}\right
\Vert +\left( 1+\frac{32M^{3}a^{3}}{\sqrt{6\pi }\left( 1-4Ma\right) }\right)
M \\ 
\left \Vert \mathfrak{u-}\underset{j=0}{\overset{n}{\sum }}\mathfrak{u}%
_{j}\right \Vert \leq \frac{16M^{2}a}{\sqrt{3\pi }\left( 1-4Ma\right) }\frac{%
\left( 4Ma\right) ^{n}}{\left( n+1\right) \sqrt{n+1}},\text{ }n\in 
\mathbb{N}
^{\ast }%
\end{array}%
\right.
\end{equation*}

\section{Proof of the main result}

Relying on the formula (\ref{fundamental formula}) we can easily prove,
under the assumption (\ref{ASSUMP1}) of the main result, that :%
\begin{equation}
||\mathfrak{u}_{n}||\leq w_{n}M^{n}a^{n-1},\text{ }n\in 
\mathbb{N}
^{\ast }  \label{Inequal}
\end{equation}%
where $\left( w_{n}\right) _{n\geq 1}$is the sequence of positive
coefficient defined by the relations :%
\begin{equation}
\left \{ 
\begin{array}{c}
w_{1}=1 \\ 
w_{n+1}=\underset{j_{1}+2j_{2}+...+nj_{n}=n}{\sum }\frac{\left(
j_{1}+...+j_{n}\right) !}{j_{1}!...j_{n}!}\left( w_{1}\right)
^{j_{1}}...\left( w_{n}\right) ^{j_{n}},\text{ }n\in 
\mathbb{N}
^{\ast }%
\end{array}%
\right.  \label{seqW}
\end{equation}%
It follows that :%
\begin{eqnarray*}
&&\underset{n\geq 1}{\sum }w_{n+1}T^{n+1}=T\underset{n\geq 1}{\sum }\left( 
\underset{j_{1}+2j_{2}+...+nj_{n}=n}{\sum }\frac{\left(
j_{1}+...+j_{n}\right) !}{j_{1}!...j_{n}!}\left( w_{1}T\right)
^{^{j_{1}}}...\left( w_{n}T^{n}\right) ^{j_{n}}\right) \\
&=&T\underset{p\geq 1}{\sum }\left( \underset{n\geq 1}{\sum }%
w_{n}T^{n}\right) ^{p}
\end{eqnarray*}%
Let us denote by $S$ the formal series $\underset{n\geq 1}{\sum }%
w_{n}T^{n}\in 
\mathbb{C}
\left[ \left[ T\right] \right] .$ Then we can write : 
\begin{equation*}
S-T=T\left( \frac{S}{1-S}\right)
\end{equation*}%
It follows that :%
\begin{equation*}
S=T\left( \frac{1}{1-S}\right)
\end{equation*}%
It follows from the well known Lagrange inversion formula for formal series (%
\cite{GESS}) that the following equalities hold for all $n\geq 1:$%
\begin{eqnarray*}
&&\left[ T^{n}\right] S=\frac{1}{n!}\left( \frac{d^{n-1}}{dt^{n-1}}\left[
\left( \frac{1}{1-t}\right) ^{n}\right] \right) _{t=0} \\
&=&\frac{1}{n!}\underset{j_{1}+...+j_{n}=n-1}{\sum }\frac{\left( n-1\right) !%
}{j_{1}!...j_{n}!}\cdot \\
&&\cdot \left( \frac{d^{j_{1}}}{dt^{j_{1}}}\left[ \left( \frac{1}{1-t}%
\right) \right] \right) _{t=0}...\left( \frac{d^{j_{n}}}{dt^{j_{n}}}\left[
\left( \frac{1}{1-t}\right) \right] \right) _{t=0} \\
&=&\frac{1}{n!}\underset{j_{1}+...+j_{n}=n-1}{\sum }\frac{\left( n-1\right) !%
}{j_{1}!...j_{n}!}j_{1}!...j_{n}!=\frac{1}{n}\underset{j_{1}+...+j_{n}=n-1}{%
\sum }1 \\
&=&\frac{\left( 2n-1\right) !}{\left( n!\right) ^{2}}
\end{eqnarray*}%
Thence we have for all $n\geq 1:$%
\begin{equation*}
w_{n}=\frac{\left( 2n-1\right) !}{\left( n!\right) ^{2}}
\end{equation*}%
The inequalities (\ref{Inequal}) become :%
\begin{equation*}
||\mathfrak{u}_{n}||\leq \frac{\left( 2n-1\right) !}{\left( n!\right) ^{2}}%
M^{n}a^{n-1},\text{ }n\in 
\mathbb{N}
^{\ast }
\end{equation*}%
But we have by virtue of the well known Stirling formula : 
\begin{equation*}
\frac{\left( 2n-1\right) !}{\left( n!\right) ^{2}}M^{n}a^{n-1}\underset{%
n\rightarrow +\infty }{\sim }\frac{1}{2\sqrt{\pi }a}\frac{1}{n\sqrt{n}}%
\left( 4Ma\right) ^{n}
\end{equation*}%
It follows then from the assumption (\ref{ASSUMP2}) that the vector series $%
\mathcal{A}dm\left( \mathfrak{N};\mathfrak{f}\right) $\ is absolutely
convergent in the Banach space $X.$ Let us then set for all $\varepsilon \in
\lbrack 0,1]:$%
\begin{equation*}
\mathfrak{U}(\varepsilon ):=\overset{+\infty }{\underset{n=0}{\sum }}%
\varepsilon ^{n}\mathfrak{u}_{n},\text{ }\mathfrak{u}:=\overset{+\infty }{%
\underset{n=0}{\sum }}\mathfrak{u}_{n}=\mathfrak{U}(1)
\end{equation*}%
Since $\mathfrak{N}$ is analytic on $E,$ it follows that we have for all $%
\varepsilon \in \lbrack 0,1]:$ 
\begin{equation*}
\mathfrak{N}(\mathfrak{U}(\varepsilon ))=\overset{+\infty }{\underset{n=0}{%
\sum }}\varepsilon ^{n}A_{n}\text{ }
\end{equation*}%
If we choose $\varepsilon =1,$ we will then obtain the relations : 
\begin{eqnarray*}
\mathfrak{N}(\mathfrak{u}) &=&\overset{+\infty }{\underset{n=0}{\sum }}A_{n}=%
\overset{+\infty }{\underset{n=0}{\sum }}\mathfrak{u}_{n+1} \\
&=&\mathfrak{u}-\mathfrak{f}
\end{eqnarray*}%
It follows that $\mathfrak{u}$ is a solution of the equation (\ref{E_CANON}%
). On the other hand we have for each $n\in 
\mathbb{N}
^{\ast }:$%
\begin{eqnarray*}
\left \Vert \mathfrak{u-}\underset{j=0}{\overset{n}{\sum }}\mathfrak{u}%
_{j}\right \Vert &\leq &\overset{+\infty }{\underset{j=n+1}{\sum }}\left
\Vert \mathfrak{u}_{j}\right \Vert \leq \overset{+\infty }{\underset{j=n+1}{%
\sum }}\frac{\left( 2j-1\right) !}{\left( j!\right) ^{2}}M^{j}a^{j-1} \\
&\leq &M\overset{+\infty }{\underset{j=n}{\sum }}\frac{\left( 2j\right)
!\left( 2j+1\right) }{\left( j!\right) ^{2}\left( j+1\right) ^{2}}\left(
Ma\right) ^{j}
\end{eqnarray*}%
But, according to (\cite{SAND}) the double inequality holds for all integers 
$j\geq 2$ :%
\begin{equation*}
\sqrt{2\pi }j^{j+\frac{1}{2}}e^{-j}\leq j!\leq \sqrt{\frac{j}{j-1}}\sqrt{%
2\pi }j^{j+\frac{1}{2}}e^{-j}
\end{equation*}%
It follows then, by easy computations, that we have for each $j\in 
\mathbb{N}
^{\ast }\backslash \left \{ 1\right \} $ :%
\begin{equation*}
\frac{\left( 2j\right) !\left( 2j+1\right) }{\left( j!\right) ^{2}\left(
j+1\right) ^{2}}\rho ^{j}\leq \frac{4}{\sqrt{3\pi }j\sqrt{j}}\left(
4Ma\right) ^{j}
\end{equation*}

\begin{itemize}
\item \bigskip First case : $4Ma=1$
\end{itemize}

In this case the following estimate holds for all $n\in 
\mathbb{N}
^{\ast }:$ 
\begin{eqnarray*}
\left \Vert \mathfrak{u-}\underset{j=0}{\overset{n}{\sum }}\mathfrak{u}%
_{j}\right \Vert &\leq &\frac{4M}{\sqrt{3\pi }}\overset{+\infty }{\underset{%
j=n+1}{\sum }}\frac{1}{j\sqrt{j}} \\
&\leq &\frac{4M}{\sqrt{3\pi }}\overset{+\infty }{\underset{n}{\int }}\frac{1%
}{s\sqrt{s}}ds \\
&\leq &\frac{8M}{\sqrt{3\pi }}\frac{1}{\sqrt{n}}
\end{eqnarray*}

It follows that :%
\begin{eqnarray*}
\left \Vert \mathfrak{u}\right \Vert &\leq &\left \Vert \mathfrak{u}%
_{0}\right \Vert +\left \Vert \mathfrak{u}_{1}\right \Vert +\left \Vert 
\mathfrak{u-}\underset{j=0}{\overset{1}{\sum }}\mathfrak{u}_{j}\right \Vert
\\
&\leq &\left \Vert \mathfrak{f}\right \Vert +M+\frac{8M}{\sqrt{3\pi }} \\
&\leq &\left \Vert \mathfrak{f}\right \Vert +\left( 1+\frac{8}{\sqrt{3\pi }}%
\right) M
\end{eqnarray*}

\begin{itemize}
\item \bigskip Second case : $4Ma<1$
\end{itemize}

In this case the following estimate holds for all $n\in 
\mathbb{N}
^{\ast }:$%
\begin{eqnarray*}
\left \Vert \mathfrak{u-}\underset{j=0}{\overset{n}{\sum }}\mathfrak{u}%
_{j}\right \Vert &\leq &\frac{4M}{\sqrt{3\pi }}\overset{+\infty }{\underset{%
j=n+1}{\sum }}\frac{1}{j\sqrt{j}}\left( 4Ma\right) ^{j} \\
&\leq &\frac{16M^{2}a}{\sqrt{3\pi }\left( 1-4Ma\right) }\frac{\left(
4Ma\right) ^{n}}{\left( n+1\right) \sqrt{n+1}}
\end{eqnarray*}%
It follows that :%
\begin{eqnarray*}
\left \Vert \mathfrak{u}\right \Vert &\leq &\left \Vert \mathfrak{u}%
_{0}\right \Vert +\left \Vert \mathfrak{u}_{1}\right \Vert +\left \Vert 
\mathfrak{u-}\underset{j=0}{\overset{1}{\sum }}\mathfrak{u}_{j}\right \Vert
\\
&\leq &\left \Vert \mathfrak{f}\right \Vert +M+\frac{64M^{4}a^{3}}{2\sqrt{%
6\pi }\left( 1-4Ma\right) } \\
&\leq &\left \Vert \mathfrak{f}\right \Vert +\left( 1+\frac{32M^{3}a^{3}}{%
\sqrt{6\pi }\left( 1-4Ma\right) }\right) M
\end{eqnarray*}

Thence the proof of the main result is complete.

$\square $

\bigskip

\bigskip

\-

\end{document}